\documentclass[oneside,english]{amsart}

\usepackage[T1]{fontenc}
\usepackage{amssymb}

\makeatletter
\numberwithin{equation}{section} %% Comment out for sequentially-numbered
\numberwithin{figure}{section} %% Comment out for sequentially-numbered
  \theoremstyle{plain}
  \newtheorem*{thm*}{Theorem}
  \theoremstyle{plain}
  \newtheorem{thm}{Theorem}[section]
  \theoremstyle{plain}
  \newtheorem{prop}[thm]{Proposition}
  \theoremstyle{remark}
  \newtheorem*{rem*}{Remark}

\usepackage{babel}
\makeatother

\begin{document}

\title{The Cauchy problem for a short-wave equation}

\author{S\'{\i}lvio M.A. Gama, Gueorgui Smirnov\\\\ Centro de Matem\'atica da Universidade do Porto \\ Departamento de
Matem\'atica Aplicada \\ Universidade do Porto\\ Portugal}



% \address{Centro de Matem\'atica da Universidade do Porto and Departamento de
% Matem\'atica Aplicada. Faculdade de Ci\^encias da Universidade do Porto.
% Rua do Campo Alegre, 687, 4169-007 Porto, Portugal.}

% \email{\{smgama,gsmirnov\}@fc.up.pt}

\begin{abstract}
We prove an existence and uniqueness of solution for the
Cauchy problem of the simplest nonlinear short-wave equation, $u_{tx}=u-3u^{2}$,
with periodic boundary condition.
\end{abstract}

\maketitle

\section{Introduction}

In this paper we consider the Cauchy problem for the short-wave equation
\begin{equation}
u_{tx}=u-3u^{2},\label{eq:sw1d}
\end{equation}
with the  periodic boundary condition ($L>0$)
\begin{equation}
u(0,t)=u(L,t),\qquad t\ge0,\label{eq:PerBoundCond}
\end{equation}
and the initial condition
\begin{equation}
u(x,0)=\phi(x),\qquad\forall x\in\mathbb{R}.\label{eq:InitCond}
\end{equation}
Here, $u(x,t)$ represents a small amplitude depending on one-dimensional (fast) space
variable $x$ and (slow) time $t$.

Nonlinear evolution of long waves in dispersive media with small amplitude
in shallow water is a well known subject and
described by many mathematical models such as
the Boussinesq equation~\cite{b,w},
the KdV equation~\cite{kdv}, or
the Benjamin-Bona-Mahony-Peregrine equation~(BBMP)~\cite{bbm,p}.
In contrast, for short-waves,
commonly called ripples, only a few results exist~\cite{mm,gkm,bkmp}.
When we speak of long or short-waves, we are referring to an underlying spacescale, $X$, to
which all space variables have been compared. Thus, for instance, for the surface-wave motion
of a fluid, the unperturbed depth serves as a natural parameter. The shortness of the waves
is referred to this underlying parameter.

The short-wave equation~\eqref{eq:sw1d} is derived in~\cite{mm}
via multiple-scale perturbation theory from BBMP and governs the leading order
term of the asymptotic
dynamics of short-waves sustained by BBMP. 
A first study of equation~(\ref{eq:sw1d}) is done in~\cite{gkm}.
Briefly, we sketch here its derivation.
Start from BBMP
\begin{equation}\label{bbmp_eq}
  U_T+U_X-U_{XXT}=3(U^2)_X,
\end{equation}
the model equation for the unperturbed
equation to which we will find the short-wave limit.
Here, $U(X,T)$ represents a small amplitude depending on one dimensional space variable $X$ and  time $T$.
Its linear dispersion relation, $\omega(k),$ is real (this means that we are not
dealing with dissipative effects) and is given by
\begin{equation}
\omega(k)=\frac{k}{1+k^2},
\end{equation}
having zero limit when $k\rightarrow \infty.$
The phase and group velocity are all bounded in the
short-wave limit $k\rightarrow \infty.$
This property allows BBMP to sustain short-waves.
In fact, let us consider a short-wave with characteristic lenght $\ell=\varepsilon\sim k^{-1},$
with $k\gg 1.$ Define the scaled (fast) space variable $x=\varepsilon^{-1}X$ ($\varepsilon\ll 1$).
The characteristic time associated with short-waves is given by looking at the dispersive
relation of the linear part for the time variable. In our case,
$\omega(\varepsilon^{-1})=\varepsilon-\varepsilon^3+\varepsilon^5-\cdots$. In this way, we obtain
the scaled (slow) time variable $t=\varepsilon T.$
We are lead thus to the scaled variables $x=\varepsilon^{-1}X$
and $t=\varepsilon T,$  which transforms the $X$ and $T$ derivatives into
$\partial_X=\varepsilon^{-1}\partial_x$ and $\partial_T=\varepsilon\partial_t.$ 
Assume now the expansion $U=u_0+\varepsilon u_1+\cdots.$ Passing to the $x$ and $t$
variables and integrating in $x,$ we have the lowest order in~(\ref{bbmp_eq}) in the form
\begin{equation}
  u_{0tx}=u_0-3(u_0)^2.
\end{equation}
For simplicity, writing $u_0$ as $u$, we obtain~(\ref{eq:sw1d}).

Under certain conditions, we prove, in the next section,
the existence and uniqueness of solutions for (\ref{eq:sw1d}) -- (\ref{eq:InitCond}).

\section{Main result}

Let $u=u(x,t)$ be a classical solution to the Cauchy problem, that is, a twice continuously differentiable
function satisfying (\ref{eq:sw1d}) - (\ref{eq:InitCond}). Integrating the left-hand side of \eqref{eq:sw1d}
in $x$, from $0$ to $L$, and using \eqref{eq:PerBoundCond}, we get 
\[
\frac{d}{dt}\int_{0}^{L}u_{x}(x,t)dx=\frac{d}{dt}\left(u(L,t)-u(0,t)\right)=0.
\]
 Therefore, from \eqref{eq:sw1d}, we have 
\begin{equation}
0=\frac{d}{dt}\int_{0}^{L}u_{x}(x,t)dx=\int_{0}^{L}\left(u(x,t)-3u^{2}(x,t)\right)dx.
\label{eq:OneRelation}
\end{equation}
Thus, it is natural to consider only initial conditions satisfying (\ref{eq:OneRelation}).

Note also that the $L_2$-norm of $u_x(\cdot, t)$ is a constant. Indeed, multiplying both sides of \eqref{eq:sw1d} 
by $u_{x}$ and integrating
in $x,$ from $0$ to $L$, we obtain 
\[
\frac{1}{2}\frac{d}{dt}|u_x(\cdot, t)|_2^2=\frac{d}{dt}\int_{0}^{L}
\frac{u_{x}^{2}(x,t)}{2}dx=\int_{0}^{L}\left(u(x,t)-3u^{2}(x,t)\right)u_{x}(x,t)dx
\]
 \[
=\int_{0}^{L}\frac{\partial}{\partial x}\left(\frac{u^{2}(x,t)}{2}-u^{3}(x,t)\right)dx
\]
\begin{equation}
\label{uxconst}
=\left(\frac{u^{2}(L,t)}{2}-u^{3}(L,t)\right)-\left(\frac{u^{2}(0,t)}{2}-u^{3}(0,t)\right)=0.
\end{equation}
This observation is of importance in the proof of a global existence. 

\vspace{2mm}

We will seek for solutions to problem (\ref{eq:sw1d}) - (\ref{eq:InitCond}) in a generalized sense. 
Namely, consider a formal Fourier series
\begin{equation}
u(x,t)=\sum_{n=-\infty}^{\infty}u_{n}(t)e^{i\frac{2\pi}{L}nx},\qquad u_{-n}=\overline{u_{n}},
\label{eq:Fourierexpansion}
\end{equation}
with coefficients depending on $t$. Assume that
\[
u(x,0)=\phi(x),\;\; x\in R,
\]
where $\phi$ is an $L$-periodic function. It is assumed that 
$u_{-n}=\overline{u_{n}}$ or, equivalently, $u(x,t)\in\mathbb{R}$.
Formally substituting Fourier series \eqref{eq:Fourierexpansion} in the differential equation we obtain  
a system of ordinary differential
equations 
\begin{equation}
\frac{du_{n}(t)}{dt}=-\frac{iL}{2\pi n}\left(u_{n}(t)-3\sum_{{\alpha+\beta=n\atop n\in\mathbb{Z}}}u_{\alpha}(t)
u_{\beta}(t)\right),\;\;\; n\neq0.\label{eq:SWfourier}
\end{equation}
(Denote $u_{n}(t)$ simply by $u_{n}.$) Note that, for $n=0,$ we
do not obain a differential equation for $u_{0},$ but a constraint
relating $u_{0}$ to all the others Fourier modes. Since $u_{0}$
is the real function $u$ average value over the domain of periodicity,  we obtain the equation
\begin{equation}
\label{u0}
u_{0}-3u_{0}^{2}=3\sum_{{n\in\mathbb{Z}\atop n\ne0}}\left|u_{n}\right|^{2}.
\end{equation}
It admits real solutions 
\begin{equation}
u_{0}=\frac{1}{6}\left(1\pm\sqrt{1-36\sum_{{n\in\mathbb{Z}\atop n\ne 0}}\left|u_{n}\right|^{2}}\right),
\label{eq:u0}
\end{equation}
only if $\sum_{{n\in\mathbb{Z}\atop n\ne0}}\left|u_{n}\right|^{2}\le1/36.$
For definiteness assume from now on that the sign in formula \eqref{eq:u0}
is plus, for example. The other choice is essentially the same, the
major difference being the fact that it results in waves travelling
in the opposite direction~\cite{gkm}.

Rewrite \eqref{eq:SWfourier}, in the integral form
\begin{equation}
\label{int}
u_{n}(t)=\phi_n-\frac{iL}{2\pi n}\int_{0}^{t}\left(u_{n}(s)-3\sum_{{\alpha+\beta=n\atop n\in\mathbb{Z}}}u_{\alpha}(s)
u_{\beta}(s)\right)ds,\;\;\; n\neq0,
\end{equation}
Denote by $H$ the
space of complex sequences $v=\{v_{n}\}_{n\in\mathbb{Z}}$ with
the norm 
\[
|v|=\left(|v_{0}|^{2}+\sum_{{n\in\mathbb{Z}\atop n\ne0}}n^{2}|v_{n}|^{2}\right)^{1/2}.\]
The space of $L$-periodic functions $u$ with the Fourier coefficients
from $H:$ $\{u_{n}\}_{n=-\infty}^{\infty}\in H$, we shall also denote
by $H$.  Let 
\[
\phi(x)=\sum_{n=-\infty}^{\infty}\phi_{n}e^{i\frac{2\pi}{L}nx}\in H,
\]
with $\phi_{-n}=\overline{\phi_{n}}$. We say that  a function $u\in C\left([0,\infty),H\right)$, 
$$
t\rightarrow u(t)=\sum_{n=-\infty}^{\infty}u_{n}(t)e^{i\frac{2\pi}{L}nx},\qquad u_{-n}=\overline{u_{n}},
$$
is a solution to problem (\ref{eq:sw1d}) - (\ref{eq:InitCond}),
if $\dot{u}\in L_{\infty}\left([0,\infty),H\right)$, 
and the Fourier coefficients $u_n$ satisfy  (\ref{eq:u0}), (\ref{int}), and $u_n(0)=\phi_n$, for all $n$.

\vspace{2mm}

Now we are in a position to formulate the main result of this paper.

\begin{thm*}
\label{th1}
If $\phi\in H$ satisfies 
$$\sum_{{n\in\mathbb{Z}\atop n\ne0}}n^{2}|\phi_{n}|^{2}<1/72
$$
 and 
$$
\int_{0}^{L}\left(\phi(x)-3\phi^{2}(x)\right)dx=0,
$$
then  problem (\ref{eq:sw1d}) - (\ref{eq:InitCond})  has one and only one solution. For all $t\geq 0,$ the 
Fourier series (\ref{eq:Fourierexpansion}) converges uniformely in $x$. Its sum is differentiable in $x$. 
The derivative satisfies the conditions $u_x(\cdot,t)\in L_2([0,L],R)$ and $u_x(x,\cdot)\in C([0,\infty[,R)$. 
Moreover, $u_x$ is differentiable in $t$ and (\ref{eq:sw1d}) holds.
\end{thm*}

We divide the proof in several steps. First note that the condition
\[
\int_{0}^{L}(\phi(x)-3\phi^{2}(x))dx=0,\]
implies 
\[
\phi_{0}=3|\phi_{0}|^{2}+3\sum_{{n\in\mathbb{Z}\atop n\ne0}}|\phi_{n}|^{2}.
\]
From this, we get 
\begin{equation}
\phi_{0}=\frac{1}{6}\left(1\pm\sqrt{1-36\sum_{{n\in\mathbb{Z}\atop n\ne0}}|\phi_{n}|^{2}}\right).\label{phi0}\end{equation}
Since 
$$
\sum_{{n\in\mathbb{Z}\atop n\ne0}}|\phi_{n}|^{2}\leq
\sum_{{n\in\mathbb{Z}\atop n\ne0}}n^2|\phi_{n}|^{2}
<1/72,
$$
 $\phi_0$ is well defined. 

Let $v(\cdot)\in L_{\infty}\left([0,T],H\right)$. The norm in this
space we shall denote by $\|v\|$. Define an operator $f:L_{\infty}\left([0,T],H\right)\rightarrow L_{\infty}\left([0,T],H\right)$
as follows: 
\begin{eqnarray}
 &  & f_{n}(v(\cdot))(t)=\phi_{n}-\frac{iL}{2\pi n}\int_{0}^{t}\left(v_{n}(s)-
3\sum_{k=-\infty}^{\infty}v_{k}(s)v_{n-k}(s)\right)ds,\;\;\; n\neq0,\label{fn}\\
 &  & f_{0}(v(\cdot))(t)=\frac{1}{6}\left(1+\sqrt{1-36\sum_{{n\in\mathbb{Z}\atop n\ne0}}|f_{n}(v(\cdot))(t)|^{2}}\right).\label{f0}
\end{eqnarray}

Let $M>0$. Denote by $\Phi\in L_{\infty}\left([0,T],H\right)$ the
constant function $\Phi(t)\equiv\phi$ and consider a complete metric
space \[
V_{TM}=\{v(\cdot)\in L_{\infty}\left([0,T],H\right)\mid\|v-\Phi\|\leq M\}\]
with the metric induced by $L_{\infty}\left([0,T],H\right)$. We need
the following auxiliary results.

\begin{prop}
\label{p1}If \[
\sum_{n\neq0}n^{2}|\phi_{n}|^{2}<1/72\]
 and $T$ is sufficiently small, then $f$ is well defined and is
a contractive map from $V_{TM}$ into $V_{TM}$. 
\end{prop}

\begin{proof}
Since \[
f_{n}(v)(t)-f_{n}(w)(t)=-\frac{iL}{2\pi n}\left.\int_{0}^{t}\right[(v_{n}(s)-w_{n}(s))\]
 \[
+3\left.\sum_{k=-\infty}^{\infty}((v_{k}(s)-w_{k}(s))v_{n-k}(s)+w_{k}(s)(v_{n-k}(s)-w_{n-k}(s))\right]ds,\;\;\ n\neq 0,
\]
 we have 
\[
\left.\left.\sum_{{n\in\mathbb{Z}\atop n\ne0}}n^{2}|f_{n}(v)(t)-f_{n}(w)(t)|^{2}\leq({\rm const})\sum_{{n\in\mathbb{Z}\atop n\ne0}}\right[\int_{0}^{t}\right[|v_{n}(s)-w_{n}(s)|\]
 \[
+\left.\left.\sum_{k=-\infty}^{\infty}(|v_{k}(s)-w_{k}(s)||v_{n-k}(s)|+|w_{k}(s)||v_{n-k}(s)-w_{n-k}(s)|\right]ds\right]^{2}\]
 \[
\leq({\rm const})t\left.\sum_{{n\in\mathbb{Z}\atop n\ne0}}\int_{0}^{t}\right[|v_{n}(s)-w_{n}(s)|\]
 \[
+\left.\sum_{k=-\infty}^{\infty}|v_{k}(s)-w_{k}(s)|(|v_{n-k}(s)|+|w_{n-k}(s)|)\right]^{2}ds\]
 \[
\leq({\rm const})t\left.\sum_{{n\in\mathbb{Z}\atop n\ne0}}\int_{0}^{t}\right[|v_{n}(s)-w_{n}(s)|^{2}\]
 \[
+\left.\left(\sum_{k=-\infty}^{\infty}|v_{k}(s)-w_{k}(s)|(|v_{n-k}(s)|+|w_{n-k}(s)|)\right)^{2}\right]ds\]
 \[
\leq({\rm const})t\left.\sum_{{n\in\mathbb{Z}\atop n\ne0}}\int_{0}^{t}\right[|v_{n}(s)-w_{n}(s)|^{2}+|v_{0}(s)-w_{0}(s)|^{2}(|v_{n}(s)|^{2}+|w_{n}(s)|^{2})\]
 \[
+\left.\left(\sum_{k\neq0}\frac{1}{k^{2}}\right)\sum_{k=-\infty}^{\infty}k^{2}|v_{k}(s)-w_{k}(s)|^{2}(|v_{n-k}(s)|^{2}+|w_{n-k}(s)|^{2})\right]ds\]
 \[
\leq({\rm const})t\int_{0}^{t}\left[\sum_{{n\in\mathbb{Z}\atop n\ne0}}|v_{n}(s)-w_{n}(s)|^{2}+|v_{0}(s)-w_{0}(s)|^{2}\sum_{{n\in\mathbb{Z}\atop n\ne0}}(|v_{n}(s)|^{2}+|w_{n}(s)|^{2})\right.\]
 \[
+\left.\sum_{k=-\infty}^{\infty}k^{2}|v_{k}(s)-w_{k}(s)|^{2}\sum_{{n\in\mathbb{Z}\atop n\ne0}}(|v_{n}(s)|^{2}+|w_{n}(s)|^{2})\right]ds\]
 \[
\leq({\rm const})t\int_{0}^{t}\left[1+\sum_{{n\in\mathbb{Z}\atop n\ne0}}(|v_{n}(s)|^{2}+|w_{n}(s)|^{2})\right]ds\|v-w\|^{2}\]
 \[
\leq({\rm const})T^{2}(1+\|v\|^{2}+\|w\|^{2})\|v-w\|^{2}.\]
 We have thus proved the following inequality \begin{equation}
\sum_{{n\in\mathbb{Z}\atop n\ne0}}n^{2}|f_{n}(v)(t)-f_{n}(w)(t)|^{2}\leq({\rm const})T^{2}(1+\|v\|^{2}+\|w\|^{2})\|v-w\|^{2}.\label{3}\end{equation}
 We also have \[
|f_{0}(v)(t)-f_{0}(w)(t)|^{2}=\frac{1}{36}\left|\sqrt{1-36\sum_{{n\in\mathbb{Z}\atop n\ne0}}|f_{n}(v)(t)|^{2}}-\sqrt{1-36\sum_{{n\in\mathbb{Z}\atop n\ne0}}|f_{n}(w)(t)|^{2}}\right|^{2}\]
$$
\leq\frac{({\rm const})\sum_{{n\in\mathbb{Z}\atop n\ne0}}(|f_{n}(v)(t)|^{2}+|f_{n}(w)(t)|^{2})}{\left|\sqrt{1-36\sum_{{n\in\mathbb{Z}\atop n\ne0}}|f_{n}(v)(t)|^{2}}+\sqrt{1-36\sum_{{n\in\mathbb{Z}\atop n\ne0}}|f_{n}(w)(t)|^{2}}\right|^{2}}
$$
\begin{equation}
\times\sum_{{n\in\mathbb{Z}\atop n\ne0}}|f_{n}(v)(t)-f_{n}(w)(t)|^{2}.
\label{4}
\end{equation}

The inclusion $v\in V_{TM}$ implies $\|v\|^{2}\leq(\|\Phi\|+\|\Phi-v\|)^{2}\leq(\|\Phi\|+M)^{2}$.
Since $\Phi=f(0)$, from (\ref{3}) we get \[
\sum_{{n\in\mathbb{Z}\atop n\ne0}}n^{2}|f_{n}(v)(t)-\phi_{n}|^{2}\leq({\rm const})T^{2}(1+(\|\Phi\|+M)^{2})^{2}.\]
 Therefore \[
\sum_{{n\in\mathbb{Z}\atop n\ne0}}|f_{n}(v)(t)|^{2}\leq2\sum_{{n\in\mathbb{Z}\atop n\ne0}}n^{2}|\phi_{n}|^{2}+2\sum_{{n\in\mathbb{Z}\atop n\ne0}}n^{2}|f_{n}(v)(t)-\phi_{n}|^{2}\]
 \[
\leq2\sum_{{n\in\mathbb{Z}\atop n\ne0}}n^{2}|\phi_{n}|^{2}+({\rm const})T^{2}(1+(\|\Phi\|+M)^{2})^{2}\leq\sigma<\frac{1}{36},\]
 whenever $T>0$ is small enough. Thus the map $f$ is well defined
(see (\ref{fn}) and (\ref{f0})). From (\ref{3}) and (\ref{4})
we obtain \[
|f_{0}(v)(t)-f_{0}(w)(t)|^{2}\leq({\rm const})\sum_{{n\in\mathbb{Z}\atop n\ne0}}|f_{n}(v)(t)-f_{n}(w)(t)|^{2}\]
 \[
\leq({\rm const})\sum_{{n\in\mathbb{Z}\atop n\ne0}}n^{2}|f_{n}(v)(t)-f_{n}(w)(t)|^{2}\leq({\rm const})T^{2}(1+\|v\|^{2}+\|w\|^{2})\|v-w\|^{2}.\]
 Invoking again (\ref{3}), we get \[
\|f(v)-f(w)\|^{2}\leq({\rm const})T^{2}(1+\|v\|^{2}+\|w\|^{2})\|v-w\|^{2}\]
 \begin{equation}
\leq({\rm const})T^{2}(1+(\|\Phi\|+M)^{2})\|v-w\|^{2}.\label{6}\end{equation}
 Specifically we have \[
\|f(v)-\Phi\|^{2}\leq({\rm const})T^{2}(1+(\|\phi\|+M)^{2})^{2}\leq M^{2},\]
 for small $T>0$. Thus we see that $f:V_{TM}\rightarrow V_{TM}$
and from (\ref{6}) it follows that $f$ is a contraction, whenever
$T>0$ is small enough. 

\end{proof}

\begin{prop}
\label{p2}Let $u\in L_{\infty}\left([0,T],H\right)$ be a solution
to the equation $u=f(u)$. Assume that \[
\sum_{{n\in\mathbb{Z}\atop n\ne0}}n^{2}|u_{n}(t)|^{2}\leq\delta<1/36.\]
 Then $u\in C\left([0,T],H\right)$ and $\dot{u}\in L_{\infty}\left([0,T],H\right)$. 
\end{prop}

\begin{proof}
Similarly to inequality (\ref{4}) we get 
\[
|u(t_{2})-u(t_{1})|^{2}=|u_{0}(t_{2})-u_{0}(t_{1})|^{2}+\sum_{{n\in\mathbb{Z}\atop n\ne0}}n^{2}|u_{n}(t_{2})-u_{n}(t_{1})|^{2}
\]
\[
\leq({\rm const})\sum_{{n\in\mathbb{Z}\atop n\ne0}}n^{2}|u_{n}(t_{2})-u_{n}(t_{1})|^{2}.
\]
 From (\ref{fn}) we see that the right side of the inequality is
less than or equal to 
\[
({\rm const})|t_{2}-t_{1}|\sum_{{n\in\mathbb{Z}\atop n\ne0}}\left|\int_{t_{1}}^{t_{2}}\left(|u_{n}(s)+3\sum_{k=-\infty}^{\infty}|u_{k}(s)||u_{n-k}(s)|\right)^{2}ds\right|
\]
 \[
\leq({\rm const})|t_{2}-t_{1}|\sum_{{n\in\mathbb{Z}\atop n\ne0}}\int_{t_{1}}^{t_{2}}\left(1+|u_{0}(s)|^{2}+\left(\sum_{{k\in\mathbb{Z}\atop k\ne0}}\frac{1}{k^{2}}\right)\sum_{{k\in\mathbb{Z}\atop k\ne0}}k^{2}|u_{k}|^{2}\right)\sum_{{n\in\mathbb{Z}\atop n\ne0}}|u_{n}(s)|^{2}ds\]
 \[
\leq({\rm const})|t_{2}-t_{1}|^{2}.\]
 This proves the continuity of $u(t)$.
Since \[
|\dot{u}_{0}|^{2}=\frac{9\left|\sum_{{n\in\mathbb{Z}\atop n\ne0}}(\dot{u}_{n}u_{-n}+u_{n}\dot{u}_{-n})\right|^{2}}{1-36\sum_{{n\in\mathbb{Z}\atop n\ne0}}|u_{n}|^{2}}\leq({\rm const})\sum_{{n\in\mathbb{Z}\atop n\ne0}}|\dot{u}_{n}|^{2}\sum_{{n\in\mathbb{Z}\atop n\ne0}}|u_{n}|^{2}\]
 and \[
\sum_{{n\in\mathbb{Z}\atop n\ne0}}n^{2}|\dot{u}_{n}|^{2}=\sum_{{n\in\mathbb{Z}\atop n\ne0}}\left(\frac{L}{2\pi}\right)^{2}\left|u_{n}-3\sum_{{n\in\mathbb{Z}\atop n\ne0}}u_{k}u_{n-k}\right|^{2},\]
 we have \[
|\dot{u}_{0}|^{2}+\sum_{{n\in\mathbb{Z}\atop n\ne0}}n^{2}|\dot{u}_{n}|^{2}\leq({\rm const})\sum_{{n\in\mathbb{Z}\atop n\ne0}}\left|u_{n}-3\sum_{{n\in\mathbb{Z}\atop n\ne0}}u_{k}u_{n-k}\right|^{2}\]
 \[
\leq({\rm const})\sum_{{n\in\mathbb{Z}\atop n\ne0}}\left(|u_{n}|^{2}+|u_{0}|^{2}|u_{n}|^{2}+\left(\sum_{{k\in\mathbb{Z}\atop k\ne0}}\frac{1}{k^{2}}\right)\sum_{{k\in\mathbb{Z}\atop k\ne0}}k^{2}|u_{k}|^{2}|u_{n-k}|^{2}\right)\]
 \[
\leq({\rm const})\left(1+|u_{0}|^{2}+\left(\sum_{{k\in\mathbb{Z}\atop k\ne0}}\frac{1}{k^{2}}\right)\sum_{{k\in\mathbb{Z}\atop k\ne0}}k^{2}|u_{k}|^{2}\right)\sum_{{n\in\mathbb{Z}\atop n\ne0}}|u_{n}|^{2}\leq({\rm const}).\]
 Thus $\dot{u}\in L_{\infty}([0,T],H)$. 
\end{proof}

\begin{rem*}
Note that we also proved that the function $u\in C\left([0,T],H\right)$ is Lipschitzian. 
\end{rem*}

Now show that generalized solutions also satisfy property (\ref{uxconst}).

\begin{prop}
\label{p3}Assume that $u\in L_{\infty}\left([0,T],H\right)$ satisfies
\eqref{eq:u0}. Then \[
\sum_{{n\in\mathbb{Z}\atop n\ne0}}n^{2}|u_{n}(t)|^{2}=({\rm const}).\]
 
\end{prop}
\begin{proof}
Indeed, we have
$$
\frac{d}{dt}
\sum_{n=-\infty}^\infty
\left(\frac{2\pi n}{L}\right)^2
|u_n|^2
=
\sum_{n=-\infty}^\infty
\left(\frac{2\pi n}{L}\right)^2
(\dot{u}_n u_{-n} + u_n\dot{u}_{-n})
$$
$$
= 
\frac{2\pi i}{L}
\sum_{n=-\infty}^\infty
n \left[ 
u_n\left( u_{-n}-3
\sum_{k=-\infty}^\infty
u_ku_{-n-k}\right) 
-
u_{-n}\left( u_{n}-3
\sum_{k=-\infty}^\infty
u_ku_{n-k}\right)
\right]
$$
$$
= -
\frac{6\pi i}{L}S,
$$
where
$$
S=
\sum_{n=-\infty}^\infty
n \left[ 
u_n
\sum_{k=-\infty}^\infty
u_ku_{-n-k} 
-
u_{-n} 
\sum_{k=-\infty}^\infty
u_ku_{n-k}
\right]. 
$$
Observe that
$$
S=
\sum_{n,k=-\infty}^\infty
n u_n u_k u_{-n-k}
-
\sum_{n,k=-\infty}^\infty
n u_{-n} u_k u_{n-k}
=
2
\sum_{n,k=-\infty}^\infty
n u_n u_k u_{-n-k}.
$$
On the other hand, introducing a new summation index $m=n-k$, we can rewrite $S$ in the following form:
$$
S=
\sum_{n,k=-\infty}^\infty
n u_n u_k u_{-n-k}
-
\sum_{m,k=-\infty}^\infty
(m+k) u_{-m-k} u_k u_{m}
=
-
\sum_{m,k=-\infty}^\infty
k u_{-m-k} u_k u_{m}.
$$
Combining this with the previous equality, we get $S=-S/2$. Thus $S=0$. 
\end{proof}

\noindent {\em Proof of Theorem.} 
From Proposition \ref{p1} we see that the problem under consideration
has one and only one solution $u\in L_{\infty}\left([0,T],H\right)$,
whenever $T>0$ is small enough. By Proposition \ref{p2} $u\in C\left([0,T],H\right)$ and $\dot{u}\in L_{\infty}\left([0,T],H\right)$. Finally, Proposition \ref{p3} implies the existence of the solution for all $t\geq0$.

Show that, $u(x,t)$, the sum of Fourier series (\ref{eq:Fourierexpansion}) satisfies (\ref{eq:sw1d}). From the inequality
$$
\sum_{{n\in\mathbb{Z}\atop n\ne0}}|u_n(t)|\leq
\sqrt{\left(\sum_{{n\in\mathbb{Z}\atop n\ne0}}\frac{1}{n^2}\right)
\sum_{{n\in\mathbb{Z}\atop n\ne0}}n^2|u_n(t)|^2}=({\rm const})
$$
we see that Fourier series (\ref{eq:Fourierexpansion}) converges uniformely in $x$ for all $t\geq 0$. The inequality
$$
\sum_{n=-\infty}^{\infty}\left|
\sum_{k=-\infty}^{\infty}u_k(t)u_{n-k}(t)\right|\leq
\sum_{k=-\infty}^{\infty}|u_k(t)|
\sum_{n=-\infty}^{\infty}|u_n(t)|
$$
implies that the series
$$
\sum_{n=-\infty}^{\infty}\left(\sum_{k=-\infty}^{\infty}u_k(t)u_{n-k}(t)\right)e^{i\frac{2\pi}{L}nx}
$$
converges for all $t\geq 0$. Multiplying (\ref{int}) by $e^
{i\frac{2\pi}{L}nx}$ and adding the obtained equalities, we get
$$
\sum_{n=-\infty}^{\infty}i\frac{2\pi}{L}n u_{n}(t)e^{i\frac{2\pi}{L}nx}
=
\sum_{n=-\infty}^{\infty}i\frac{2\pi}{L}n \phi_n e^{i\frac{2\pi}{L}nx}
$$
$$
+
\sum_{{n\in\mathbb{Z}\atop n\ne0}}
\int_{0}^{t}\left(u_{n}(s)-3\sum_{{\alpha+\beta=n\atop n\in\mathbb{Z}}}u_{\alpha}(s)u_{\beta}(s)\right) e^{i\frac{2\pi}{L}nx}ds
$$
From the Lebesgue theorem and the above estimates we have
$$
u_x(x,t)
=
\phi_x(x)
+\int_{0}^{t}
\sum_{{n\in\mathbb{Z}\atop n\ne0}}
\left(u_{n}(s)-3\sum_{{\alpha+\beta=n\atop n\in\mathbb{Z}}}u_{\alpha}(s)u_{\beta}(s)\right) e^{i\frac{2\pi}{L}nx}ds.
$$
Combining this with (\ref{u0}) we obtain
$$
u_x(x,t)
=
\phi_x(x)
+\int_{0}^{t}
(u(x,s)-3u^2(x,s))ds.
$$
This ends the proof. \hfill $\Box$


\begin{thebibliography}{1}

\bibitem{bbm}
T.B. Benjamin, J.L. Bona and J.J. Mahony,
\textit{Model Equations for Long Waves in Nonlinear Dispersive Systems},
Phil. Trans. R. Soc. A \textbf{272} (1972), 47 -- 78.


\bibitem{bkmp}
C.H. Borzi, R.A. Kraenkel, M.A. Manna and A. Pereira,
\textit{Nonlinear dynamics of short travelling capillary-gravity waves},
Phys. Rev. E \textbf{71}(2) (2005), 026307-1 -- 026307-9.


\bibitem{b}
J. Boussinesq,
\textit{Th\'eorie de l'intumescence liquide, applel\'ee onde
solitaire ou de translation, se propageant dans un canal rectangulaire},
Compte Rendue Acad. Sci. Paris \textbf{72} (1871), 755 -- 759.


\bibitem{gkm}
S.M. Gama, R.A. Kraenkel and M.A. Manna,
\textit{Short-waves instabilities in the Benjamin-Bona-Mahoney-Perigrine
equation: theory and numerics},
Inverse Problems \textbf{17}(4) (2001), 864 -- 870.


\bibitem{kdv}
D.J. Korteweg and G. deVries,
\textit{On the Change of Form of Long Waves advancing
in a Rectangular Canal and on a New Type of
Long Stationary Waves},
Philos. Mag. \textbf{36}(5) (1895), 422 -- 443.


\bibitem{mm}
M.A. Manna and V. Merle
\textit{Asymptotic dynamics of short waves in nonlinear
dispersive models},
Phys. Rev. E \textbf{57}(5) (1998), 6206 -- 6209.


\bibitem{p}
D.H. Peregrine,
\textit{Long waves on a beach},
J. Fluid Mech. \textbf{27} (1967), 815 -- 827.


\bibitem{w}
G.B. Whitham
\textit{Linear and Nonlinear Waves}, (1972),
Wiley Interscience, New York.



\end{thebibliography}
\end{document}